\documentclass[13pt,twoside]{article}
\usepackage{amsmath, amsthm, amsfonts, amssymb, color, mathrsfs, dsfont}
\usepackage{mathrsfs}

\def\az{\alpha}  \def\bz{\beta}
    \def\dz{\delta}
    \def\fz{\varphi}
  
\def\lz{\lambda} \def\mz{\mu}
\def\nz{\nu}     
\def\pz{\pi}

        \def\sz{\sigma}
        
\def\vz{\varepsilon}

\def\llz{\Lambda}

\def\qd{\quad}
\def\qqd{\qquad}

\def\lt{\left}
\def\rt{\right}

\newcommand{\mathsym}[1]{{}}

\def\leq{\leqslant}
\def\geq{\geqslant}

\def\nnd{\noindent}

\def\thm{\nnd\bg{thm1}}

\def\lmm{\nnd\bg{lmm1}}
\def\prp{\nnd\bg{prp1}}

\def\xmp{\nnd\bg{xmp1}}

\def\defn{\nnd\bg{defn1}}

\def\dethm{\end{thm1}}

\def\delmm{\end{lmm1}}
\def\deprp{\end{prp1}}
\def\dexmp{\end{xmp1}}

\def\dedefn{\end{defn1}}

\def\prf{\medskip \noindent {\bf Proof}. }
\def\deprf{\quad $\square$ \medskip}
\def\bg{\begin}
\def\be{\bg{equation}}
\def\de{\end{equation}}

\def\dear{\end{eqnarray}}
\def\lb{\label}

\def\den{\end{enumerate}}

\def\var{\text{\rm Var}}

\def\d{\text{\rm d}}
\def\gap{\text{\rm gap}}

\begin{document}

\allowdisplaybreaks[4]
\thispagestyle{empty}
\renewcommand{\thefootnote}{\fnsymbol{footnote}}


\vspace*{.5in}
\begin{center}
{\bf\Large Estimate the exponential convergence rate of $f$-ergodicity via spectral gap } \\
{
Xianping Guo\footnote{School of Mathematics, Sun Yat-Sen University, China. E-mail: mcsgxp@mail.sysu.edu.cn}
Zhong-Wei Liao\footnote{South China Research Center for Applied Mathematics and Interdisciplinary Studies, South China Normal University, China. E-mail: zhwliao@hotmail.com},
}
\end{center}


\bigskip

\noindent {\bf Abstract} \qd This paper studies the $f$-ergodicity and its exponential convergence rate for continuous-time Markov chain. Assume $f$ is square integrable, for reversible Markov chain, it is proved that the exponential convergence of $f$-ergodicity holds if and only if the spectral gap of the generator is positive. Moreover, the convergence rate is equal to the spectral gap. For irreversible case, the positivity of spectral gap remains a sufficient condition of $f$-ergodicity. The effectiveness of these results are illustrated by some typical examples.


\vskip.1in


\medskip

\section{Introduction and main results}

In this paper, we study the rate of convergence to equilibrium of continuous-time Markov chain. Assume $(X_t)_{t \geq 0}$ is a positive recurrent Markov chain defined on a countable state space $E$ with stationary distribution $\pz$. Denote by $Q= (q_{ij})$ and $P_t (i, j)$ the $Q$-matrix and the corresponding Markov semigroup. For any measurable function $f: E \rightarrow [1, \infty)$, the $f$-norm of signed measure $\mz$ is defined as $\| \mz\|_f := \sup_{|g|\leq f} |\mz (g)|$. When $f$ is a constant function, the $f$-norm is nothing but the total variation norm. The main objective is the $f$-ergodicity of $P_t$, which means that for all $i \in E$, we have
\be\lb{def f-erg}
\lim_{t\rightarrow \infty} r (t) \| P_t (i, \cdot) - \pz \|_f =0,
\de
where $f$ satisfies $\pz(f)< \infty$ and $r(t)$ is a positive function being used to describe the convergence rate. For example, the exponential convergence means $r(t) = e^{\vz t}$, $\vz > 0$. Refer to \cite[Chapter 14]{Meyn} or \cite{Tuominen} for more details about the terminology and notations.

For $f\equiv 1$, (\ref{def f-erg}) depicts the long time behavior of Markov semigroup in total variation norm. There are many approaches in the quantitative research, refer to \cite{Chen1998}, \cite{Chen2000}, \cite{Down} and \cite{Lund}. For example, one of the efficient instrument popularized by Meyn and Tweedie is the drift condition (or Foster-Lyapunov control conditions), which implies the exponential convergence, see \cite{Down} or \cite{Meyn}. Another useful tool is functional inequalities. Assume the semigroup $P_t$ is reversible with respect to $\pz$, which means $\pz_i P_t (i, j) = \pz_j P_t (j, i)$, for all $i, j \in E$ and $t \geq 0$ (equivalently, $ \pz_i q_{ij} = \pz_j q_{ji}$). The Poincar\'{e} inequality is defined as
\be\lb{PI}
C_{\rm PI} \var_{\pz} (g) \leq (- Q g, g), \qqd g \in L^2 (\pz),
\de
where $C_{\rm PI}$ is denotes as the optimal constant and $(\cdot, \cdot)$ is the inner product in $L^2 (\pz)$. This inequality is also referred to spectral gap inequality, since the spectral gap of $Q$ can be redefined as the optimal constant of the Poincar\'{e} inequality:
$$
\gap (Q) = \inf \{ (-Q g, g): \text{$\pz(g) =0$ and $\| g \|_{L^2(\pz)} =1$} \}.
$$
Corresponding to the spectral gap is the exponential ergodicity in $L^2 (\pz)$:
$$
\lt\| P_t g - \pz(g) \rt\|_{L^2 (\pz)} \leq e^{- C_{\rm PI} t} \lt\| g- \pz(g) \rt\|_{L^2 (\pz)}, \qqd g\in L^2(\pz).
$$
According to Cauchy-Schwarz inequality, the Poincar\'{e} inequality implies the exponential convergence in total variation distance. There is a great deal of publications in this field, see for instance \cite{Bakry14}, \cite{Chen2004}, \cite{Liggett} and references within. In addition, the relationship between the Mayn-Tweedie approach and the functional inequality approach has been discussed in \cite{Bakry07}. Generalizations of functional inequalities have been studied by several authors, here we refer to \cite{Cattiaux}, \cite{Mao} and \cite{Rockner} for related results on weak Poincar\'{e} inequalities and weak logarithmic Sobolev inequalities.

For $f\geq 1$ and $r \equiv 1$ in (\ref{def f-erg}), that is the $f$-ergodicity introduced in \cite[Chapter 14]{Meyn}, but without consideration the convergence rate. What we concern is the case $f\geq 1$ and $r(t) = e^{-\vz t}$, in other words, the semigroup $P_t$ is said to have exponential $f$-ergodicity if there exists constants $\vz>0$ and $C(i, f) >0$ such that
\be\lb{f-erg}
\| P_t (i, \cdot) - \pz \|_f \leq C(i, f) e^{- \vz t}, \qqd \forall i \in E, \ t \geq 0.
\de
The maximal parameter $\vz_{\rm max}$ is called the exponential convergence rate of $f$-ergodicity. Our objectives in this paper are the criterion of $f$-ergodicity and the estimation of the convergence rate in (\ref{f-erg}) .

Researches surrounding $f$-ergodicity is applied in the theory of controlled Markov models (Markov decision processes) in \cite{Guo}. Specifically, it ensures the existence of average optimal policies in the unbounded rewards model. Hence, the explicit criterion of $f$-ergodicity is the original motive of this thesis. Following the Meyn-Tweedie approach, Douc et al. \cite{Douc} give a general form of drift condition, which is depend on the notion of ``petite set".

The main tool we use is the functional inequality. We review the conditions of $f$ in (\ref{def f-erg}). The condition ``$f\geq 1$" ensures that the $f$-ergodicity of $P_t$ implies the original ergodicity. However, it is not essential because it can be replaced by ``$f\geq \dz$" for any $\dz >0$. The condition ``$\pz(f) < \infty$ " is necessary, otherwise (\ref{def f-erg}) might be not well-define. Furthermore, when $f \in L^2 (\pz)$, the exponential rate of $f$-ergodicity $\vz_{\rm max}$ can be estimated by the spectral gap of generator, which is our main result. Different from the drift conditions given in \cite{Lund} and \cite{Meyn}, we introduce a new equivalent condition of $f$-ergodicity. The principal tools are Poincar\'{e} inequality and $h$-transform, which will be given in Section 3.
\thm\lb{Thm f-ergodic1}
{
Assume $f \in L^2 (\pz)$ and $P_t$ is reversible. Then $P_t$ has exponential $f$-ergodicity if and only if the spectral gap of $Q$-matrix $\gap (Q)>0$. Moreover, the convergence rate satisfies $\vz_{\rm max} = \gap (Q)$, and constant of (\ref{f-erg}) is $C(i, f) = \pz \lt( f^2 \rt)^{1/2} \lt( \pz^{-1}_i  -1 \rt)^{1/2}$.
}
\dethm
Since $(X_t)_{t \geq 0}$ is positive recurrent, the stationary distribution $\pz$ satisfies $\pz_i > 0$, $\forall i \in E$. Hence $C(i, f) < \infty$, $\forall i \in E$. For irreversible case, the above-mentioned equivalence will be false. However, the spectral gap condition is still a sufficient condition of $f$-ergodicity.
\prp\lb{crl f-ergodic}
{
Assume $\pz (f^2) < \infty$. If the semigroup $P_t$ is irreversible, then $\gap (Q)>0$ implies the exponential $f$-ergodicity of $P_t$.
}
\deprp

In Section 2, we will give some examples to illustrate the effectiveness of Theorem \ref{Thm f-ergodic1}. The $h$-transformation will introduced in Section 3, and then we give the proof of Theorem \ref{Thm f-ergodic1} by this method.

\section{Examples}

As previous mentioned, one of the practical criterion is the drift condition (cf. \cite{Lund}). In practical applications, this criterion is easy to verify, although the invariant measure is unknown. However, the next example show that the drift condition can not give a exact estimation of the convergence rate $\vz_{\rm max}$.

\xmp
{\rm
 Take$E= \mathds{Z}_+$. Let $\pz_i >0, (\forall i \in E)$ be an arbitrary distribution on $E$ and the $Q$-matrix defined as: $q_{ij}= \pz_{j}$ for $j\neq i$; $q_{ii}= - \sum_{j\neq i} q_{ij}$. Consider the $f$-ergodicity with $f$ satisfying $f_0 =1$, $f_i \equiv \bz > 1$, $\forall i \geq 1$. Theorem \ref{Thm f-ergodic1} gives the convergence rate as
$$
\vz_{\rm max} =1 \qd \text{and} \qd C(i, f)= \lt[ \pz_0 + \bz^2 (1- \pz_0) \rt]^{1/2} \lt(\frac{1}{\pz_i} -1 \rt)^{1/2}.
$$
}
\dexmp

\prf By the definition of $Q$-matrix, we have $\pz Q =0$ and $\pz_{i} q_{ij}= \pz_{j} q_{ji}$, which means $Q$-matrix is reversible respect to the stationary distribution $\pz$. For any $g \in L^2(\pz)$ we have
$$
(- Q (g), g) = \sum_j \pz_j g_j \lt(\sum_k q_{jk} (g_k -g_j) \rt) = \pz(g)^2 - \pz(g^2) = \var_{\pz} (g).
$$
Hence $\gap(Q)=1$. Since $f\in L^2 (\pz)$, using Theorem \ref{Thm f-ergodic1}, we obtain
$$
\vz_{\rm max} =1, \qqd C(i, f)= \lt[ \pz_0 + \bz^2 (1- \pz_0) \rt]^{1/2} \lt(\frac{1}{\pz_i} -1 \rt)^{1/2}.
$$
This gives the exact description of $\vz_{\rm max}$ and $C (i, f)$.

If we use the drift condition given in \cite[Theorem 2.2]{Lund}, we need to solve the equation
$$
Q f (i) \leq -c f (i) + b \mathds{1}_{\{0\}} (i), \qqd \forall i \in E.
$$
For $x \neq 0$, it implies that $c \leq 1- \bz^{-1} \pz (f) = \pz_0 (1- \bz^{-1}) < 1$. For $x= 0$, we have $b\geq \bz(1- \pz_0) + (\pz_0+ c -1)$. Hence the drift condition shows the $f$-ergodicity holds with convergence rate $c\in (0, \pz_0(1- \bz^{-1})]$.
\deprf

For the irreversible case, the equivalence in Theorem \ref{Thm f-ergodic1} is erroneous. The following example show that there is some difference between $\gap(Q)$ and $\vz_{\rm max}$.

\xmp {\bf (irreversible case)}
{\rm
Let $E=\{0, 1, 2 \}$ and $f$ satisfying $f_i \in [1, \infty)$, $i=1, 2, 3$. Consider the process with $Q$-matrix
$$
Q=
\left(
  \begin{array}{ccc}
    -1/2 & 1/2 & 0 \\
    0 & -1 & 1 \\
    1 & 0 & -1 \\
  \end{array}
\right).
$$
Then the process has $f$-ergodicity with convergence rate $5/4 > \gap(Q) =1$.
}
\dexmp

\prf By $\pz Q =0$, we have $\pz_0 =1/2$, $\pz_1= \pz_2= 1/4$. In this irreversible situation, we adopt the $Q$-matrix by the symmetrizing procedure. Let $\hat{q}_{ij} = \pz_j q_{ji}/ \pz_i$ and $\bar{q}_{ij} = (q_{ij}+ \hat{q}_{ij})/2$, then we have
$$
\hat{Q}=
\left(
  \begin{array}{ccc}
    -1/2 & 0 & 1/2 \\
    1 & -1 & 0 \\
    0 & 1 & -1 \\
  \end{array}
\right), \qqd
\bar{Q}=
\left(
  \begin{array}{ccc}
    -1/2 & 1/4 & 1/4 \\
    1/2 & -1 & 1/2 \\
    1/2 & 1/2 & -1 \\
  \end{array}
\right)
$$
The matrix $\bar{Q}$ is symmetry with respect to $\pz$ and it is easy to calculate that $\gap(\bar{Q}) = \gap(Q) =1>0$. Hence, by Proposition \ref{crl f-ergodic}, we obtain the $f$-ergodicity of this $Q$-process and $\vz_{\rm max} \geq 1$. However, we can not get the exact value of the convergence rate from Proposition \ref{crl f-ergodic}.

Fortunately, the convergence rate of $f$-ergodicity could be calculated directly. Firstly, the eigenvalues of $Q$ are
$$
\lz_0 =0, \qd \lz_1 = - \frac{5}{4} + \frac{\sqrt{7} i}{4}, \qd \lz_2 = - \frac{5}{4} - \frac{\sqrt{7} i}{4}.
$$
By the representation $P_t = U \llz_t U^{-1}$, where $U$ is a matrix whose column vectors are the eigenvectors, $\llz_t$ is a diagonal matrix $\llz_t= {\rm diag}({\rm exp}(\lz_i t))$, then
$$
P_t = e^{- 5/4 t} R_t +
\left(
  \begin{array}{ccc}
    1/2 & 1/4 & 1/4 \\
    1/2 & 1/4 & 1/4 \\
    1/2 & 1/4 & 1/4 \\
  \end{array}
\right).
$$
where
\begin{align*}
R_t = & \frac{\sin(\sqrt{7} t)}{\sqrt{7}}
\left(
  \begin{array}{ccc}
    -1/2 & 1/2 & 0 \\
    0 & -1 & 1 \\
    1 & 0 & -1 \\
  \end{array}
\right)  \\
& \qd + \lt( \cos(\sqrt{7} t) + \frac{5 \sin(\sqrt{7} t)}{\sqrt{7}} \rt)
\left(
  \begin{array}{ccc}
    1/2 & -1/4 & -1/4 \\
    -1/2 & 3/4 & -1/4 \\
    -1/2 & -1/4 & 3/4 \\
  \end{array}
\right)
\end{align*}
By the representation of $P_t$ and (\ref{f_norm}), we can calculate the convergence rate immediately
\begin{align*}
\|P_t (x, \cdot) - \pz \|_f = \sum_{i=0}^2 \lt| f_i \lt(p_t (x, i) - \pz_i \rt)\rt| = e^{-5/4t} \lt| \sum_{i=0}^2 f_i R_t (x, i) \rt|.
\end{align*}
The convergence rate of $f$-ergodicity is $5/4$ which is bigger than $\gap(Q)$.
\deprf

\section{The proofs}

The $h$-transform (or Doob's $h$-transform) is an useful transformation in probability or potential theory. For instance, in \cite{Pinsky}, the principal eigenvalue of diffusion operators have been carefully handled by the $h$-transform and applied to multi-dimensional case. Refer to \cite[Chapter 1]{Bakry14} for more details.

Let $P_t$ be a Markov semigroup with stationary measure $\pz$ and $f \in L^2 (\pz)$ be a strictly positive measurable function. Define a new semigroup as
$$
P_t^f (g) = \frac{1}{f} P_t (fg), \qqd \forall g \in L^2(\pz), \  t\geq 0.
$$
Similarly, the $h$-transform of $Q$-matrix and stationary distribution $\pz$ are
$$
Q^f (g) (i) = \frac{1}{f(i)} \sum_{j \in E} q_{ij} f(j) g(j), \qd \pz^f (g) (i):= \dfrac{1}{f (i)} \sum_{j \in E} \pz_j f (j) g (j), \qd \forall g \in L^2(\pz).
$$
When $P_t$ is reversible, it is easy to show that $P_t^f $ is reversible with respect to measure $\nz_i := f^2 (i) \pz_i $. Moreover, the semigroup $P_t^f$ has similar properties with $P_t$.
\lmm\lb{f-erg lmm1}
{
Let $P_t$ be a reversible Markov semigroup with respect to $\pz$, define $P_t^f$, $\pz^f$ and $\nz$ as mentioned above. For any function $g_1, g_2 \in L^2(\nz)$ we have:

(1) Semigroup property: $P_{t+ s}^f = P_t^f P_s^f , \qd \forall t, s \geq 0$;

(2) Conjugacy: $(P_t^f g_1, g_2)_{\nz}= (g_1, P_t^f g_2)_{\nz}$; $\lt(\pz^f (g_1), g_2\rt)_{\nz}= \lt(g_1, \pz^f (g_2) \rt)_{\nz}$;

(3) $\pz^f \lt(P_t^f g_1 \rt) = P_t^f \lt(\pz^f (g_1) \rt) = \pz^f (g_1)$.
}
\delmm
The proof of Lemma \ref{f-erg lmm1} is easy and straightforward. It should be noted that the $P_t^f$ is not a Markov semigroup though its properties are similar to $P_t$, and $\nz$ is not a probability measure. In order to ensure $\nz$ to be a finite measure, we need the condition $f\in L^2 (\pz)$.

This section is devoted to prove an equivalence of the exponential $f$-ergodicity and the exponential convergence of the semigroup $P_t^f$. We will start with reversible case. The irreversible case can be reduced to the symmetric one, which will be discussed shortly in the end this section.

\defn\lb{Def l2 con}
The semigroup $P_t^f$ converges exponentially in the $L^2 (\nz)$-norm if there is a constant $\sz > 0$ such that
\be\lb{L2 exp}
\lt\| P_t^f g - \pz^f (g) \rt\|_{L^2 (\nz)} \leq \lt\| g- \pz^f (g) \rt\|_{L^2 (\nz)} e^{- \sz t} , \qd \forall t\geq 0, \ g\in L^2 (\nz).
\de
The largest $\sz$ is denoted by $\sz_{\rm max}$, which is called the $L^2 (\nz)$-exponential convergence rate.
\dedefn
It is known that the exponential ergodicity rate in total variation norm (when $f\equiv 1$) is given by the spectral gap of the $Q$-matrix, refer to \cite{Chen2004}. Hence, it is natural to consider the relationship of $f$-ergodicity and the spectral gap of $Q^f$. Firstly, we give this equivalence between the convergence of $P_t^f$ and the $f$-ergodicity, which is inspired by the $h$-transform and \cite[Theorem 9.15]{Chen2004}.
\thm\lb{Thm f-ergodic}
{
Assume that $\pz (f^2) < \infty$ and $P_t$ is reversible. Then $P_t$ satisfies exponential $f$-ergodicity if and only if the semigroup $P_t^f$ converges exponentially in the $L^2 (\nz)$-norm. Moreover, we have $\vz_{\rm max} = \sz_{\rm max}$.
}
\dethm

To begin with, we give some short lemmas about the operator norm of $P_t^f$.

\lmm\lb{f-erg lmm2}
{
Let $P_t$ be a reversible semigroup. Define $P_t^f$, $\pz^f$ and $\nz$ as mentioned above, then we have
$$
\lt\| P_t^f - \pz^f \rt\|_{L^{\infty} (\nz) \rightarrow L^2 (\nz)}^2 = \lt\| P_{2t}^f - \pz^f \rt\|_{L^{\infty} (\nz) \rightarrow L^1 (\nz)}.
$$
}
\delmm

\prf For any $g \in L^{\infty} (\nz)$, by the semigroup property and conjugacy of $P_t^f$ in Lemma \ref{f-erg lmm1}, we have
\begin{align*}
\Big\| \Big( P_t^f &- \pz^f \Big) g \Big\|_{L^2 (\nz)}^2 = \lt( g, \Big( P_t^f - \pz^f \Big)^2 g \rt)_{\nz} = \lt( g, \Big( P_{2t}^f - \pz^f \Big) g \rt)_{\nz} \\
&\leq \| g \|_{L^{\infty} (\nz)} \lt\| \Big( P_{2t}^f - \pz^f \Big) g \rt\|_{L^1 (\nz)} \leq \| g \|_{L^{\infty} (\nz)}^2 \lt\| P_{2t}^f - \pz^f \rt\|_{L^{\infty} (\nz) \rightarrow L^1 (\nz)}.
\end{align*}
The last inequality gives $\lt\| P_t^f - \pz^f \rt\|_{L^{\infty} (\nz) \rightarrow L^2 (\nz)}^2 \leq \lt\| P_{2t}^f - \pz^f \rt\|_{L^{\infty} (\nz) \rightarrow L^1 (\nz)}$.

The inverse inequality is obvious by the conjugacy of $P_t^f$, details as below
\begin{align*}
\Big\| P_{2t}^f &- \pz^f \Big\|_{L^{\infty} (\nz) \rightarrow L^1 (\nz)} \leq \lt\| P_{t}^f - \pz^f \rt\|_{L^{\infty} (\nz) \rightarrow L^2 (\nz)} \lt\| P_{t}^f - \pz^f \rt\|_{L^2 (\nz) \rightarrow L^1 (\nz)} \\
&= \lt\| P_{t}^f - \pz^f \rt\|_{L^{\infty} (\nz) \rightarrow L^2 (\nz)} \lt\| \lt( P_{t}^f - \pz^f \rt)^* \rt\|_{L^{\infty} (\nz) \rightarrow L^2 (\nz)} \\
&= \lt\| P_{t}^f - \pz^f \rt\|_{L^{\infty} (\nz) \rightarrow L^2 (\nz)}^2,
\end{align*}
here $\lt(P_t^f - \pz^f \rt)^*$ is the dual of $P_t^f - \pz^f$ with respect to $(\cdot, \cdot)_{\nz}$.
\deprf

The next lemma is about the relationship between the operator norm of $P_t^f$ and the $f$-ergodicity.

\lmm\lb{f-erg lmm3}
{
Under the same conditions of Lemma \ref{f-erg lmm2}, we have
\be\lb{lmm2 =2}
\lt\| P_t^f - \pz^f \rt\|_{L^{\infty} (\nz) \rightarrow L^1 (\nz)} \leq \sum_{i \in E} \pz_i f (i) \lt\| P_{t} (i, \cdot) - \pz \rt\|_f.
\de
}
\delmm

\prf
For any $g\in L^{\infty} (\nz)$, we have $g/ \| g\|_{L^{\infty} (\nz)} \leq 1$. Directly calculating, we have
\begin{align*}
\sum_{i \in E} & \pz_i f(i) \| P_t (i, \cdot) - \pz \|_f = \sum_{i \in E} \pz_i f (i) \sup_{|\fz| \leq f} \lt| \lt( P_t - \pz \rt) (\fz) (i) \rt| \\
&= \sum_{i \in E} \pz_i f^2 (i) \sup_{|\fz|/f \leq 1} \lt| \lt( P_t^f - \pz^f \rt) \lt( \frac{\fz}{f} \rt) (i) \rt| \\
&= \sum_{i \in E} \pz_i f^2 (i) \sup_{|\fz^*| \leq 1} \lt| \lt( P_t^f - \pz^f \rt) (\fz^*) (i) \rt| \qqd \text{( where $\fz^*:= \fz/f$ )} \\
&\geq \sum_{i \in E} \pz_i f^2 (i) \lt| \lt( P_t^f - \pz^f \rt) \lt(\frac{g}{\|g \|_{L^\infty (\nz)}} \rt) (i) \rt| = \frac{ \lt\| \lt(P_t^f -\pz^f \rt) g \rt\|_{L^1 (\nz)}}{\|g \|_{L^\infty (\nz)}}.
\end{align*}
That implies
\begin{align*}
\lt\|P_t^f - \pz^f \rt\|_{L^{\infty} (\nz) \rightarrow L^1 (\nz)} &= \sup_{g \in L^{\infty} (\nz)} \frac{\lt\| \lt(P_t^f -\pz^f \rt) g \rt\|_{L^1 (\nz)}}{\| g \|_{L^{\infty} (\nz)}} \\
&\leq \sum_{i \in E} \pz_i f (i) \lt\| P_{t} (i, \cdot) - \pz \rt\|_f. \qqd \qqd \qqd  \text{\deprf}
\end{align*}

By Hahn decomposition theorem, every signed measure $\nz$ has a unique decomposition into a difference $\nz= \nz^+ - \nz^-$ of two positive measures $\nz^+$ and $\nz^-$, then the total variation norm of $\nz$ is given simply by
$$
\| \nz \|_{\rm var} = \sup_{|g|\leq 1} |\nz(g)| = \sum_{i \in E} |\nz_i|,
$$
where $|\nz| := \nz^+ + \nz^-$. Therefore, for any positive function $f$,
\be\lb{f_norm}
\| \nz \|_f = \sup_{|g|\leq f} |\nz(g)| = \sup_{|g|/f \leq 1} \lt|\nz \lt[ f \lt(\frac{g}{f} \rt) \rt] \rt| = \sum_{i \in E} f (i) |\nz_i| .
\de
Furthermore, we have following lemma.

\lmm\lb{f-erg lmm4}
{
For any probability measure $\mz$, define $h_i = \mz_i / \pz_i$. Then we have
$$
\| \mz P_t - \pz \|_f = \lt\| f \lt( P_t^* \lt( h \rt) -1  \rt)  \rt\|_{L^1 (\pz)},
$$
where $P_t^*$ is the dual semigroup of $P_t$, which means $P_t^* (i, j) := P_t (j, i) \pz_j / \pz_i$. If $P_t$ is reversible respect to $\pz$, we have $P_t^* = P_t$.
}
\delmm
\prf
The proof is straightforward. By the Hahn decomposition, we have
\begin{align*}
\| \mz P_t &- \pz \|_f = \sum_{i \in E} f (i) \lt|(\mz P_t) (i) - \pz_i \rt| = \sum_{i \in E} f(i) \lt| \sum_{j \in E} \mz_j P_t (j, i) - \pz_i \rt| \\
&= \sum_{i \in E} f(i) \lt| \sum_{j \in E} h_j \pz_j P_t (j, i)  - \pz_i \rt| = \sum_{i \in E} f(i) \lt| \sum_{j \in E} h_j \pz_i P_t^* (i, j) - \pz_i \rt| \\
&= \sum_{i \in E} \pz_i f(i) \lt| P_t^* (h) (i)  -1 \rt|. \qqd \qqd \text{\deprf}
\end{align*}

Having these preparations at hand, we are ready to prove the main results.

{\bf Proof of Theorem \ref{Thm f-ergodic}.} (i). We consider the sufficiency of Theorem \ref{Thm f-ergodic}. Assume $P_t$ satisfies exponential $f$-ergodicity, which means that there exists constants $\vz_{\rm max}>0$ and $C (i, f)>0$ such that (\ref{f-erg}) holds. Firstly, we give a direct proof under a technical condition:
\be\lb{technical}
\pz(f C (\cdot, f)) = \sum_{i \in E} \pz_i f(i) C (i, f) < \infty.
\de
By Lemma \ref{f-erg lmm2}, Lemma \ref{f-erg lmm3} and (\ref{technical}), we have
\begin{align*}
\Big\| P_t^f &- \pz^f \Big\|_{L^{\infty} (\nz) \rightarrow L^2 (\nz)}^2 = \Big\| P_{2t}^f - \pz^f \Big\|_{L^{\infty} (\nz) \rightarrow L^1 (\nz)} \\
&\leq \sum_{i \in E} \pz_i f(i) \lt\| P_{2t} (i, \cdot) - \pz \rt\|_f \leq \pz(f C (\cdot, f)) e^{-2 \vz_{\rm max} t}.
\end{align*}
Hence, for any $g$ satisfies $g\in L^{\infty} (\nz)$ and $\nz (g^2) =1$, we have
$$
\lt\| \lt( P_t^f - \pz^f \rt) g \rt\|_{L^2 (\nz)}^2 \leq \pz(f C (\cdot, f)) \lt\| g \rt\|_{L^{\infty} (\nz)}^2 e^{-2 \vz_{\rm max} t}.
$$

The constant $\pz(f C (\cdot, f))\lt\| g \rt\|_{L^{\infty} (\nz)}^2$ in the last line can be removed, which is inspired by \cite{Wang00}. For every $g$ with $\pz (f g) =0$ and $\nz (g^2) =1$, using \cite[Lemma 2.2]{Wang00} and the spectral representation theorem, we have
\begin{align*}
\lt\| P_t^f  g \rt\|_{L^2 (\nz)}^2 &= \lt\|P_t (fg) \rt\|_{L^2 (\pz)}^2 = \int_{0}^\infty e^{-2 \az t} \d \lt( E_{\az} (f g), f g \rt) \\
&\geq \lt[ \int_{0}^\infty e^{-2 \az s} \d \lt(E_{\az} (f g), f g \rt) \rt]^{t/s} \qd \text{($\forall s\leq t$, by Jensen's inequality)} \\
&= \lt\| P_s (f g) \rt\|_{L^2 (\pz)}^{2t/s} = \lt\| P_s^f g \rt\|_{L^2 (\nz)}^{2t/s},
\end{align*}
where $E_{\az}$ is the spectral measure of the generator with respect to $\az$. Thus,
$$
\lt\| P_s^f  g \rt\|_{L^2 (\nz)}^2 \leq \lt[ \pz(f C (\cdot, f)) \| g \|_{L^{\infty} (\nz)}^2 \rt]^{s/t} e^{-2 \vz_{\rm max} s}.
$$
Letting $t \rightarrow \infty$, we obtain
$$
\lt\| P_s^f  g \rt\|_{L^2 (\nz)}^2 \leq e^{-2 \vz_{\rm max} s}, \qqd \text{$g\in L^{\infty} (\nz)$, $\| g \|_{L^2 (\nz)} =1$, $\pz (f g) =0$}.
$$
Finally, since $L^{\infty} (\nz)$ is dense in $L^2 (\nz)$, we have $\vz_{\rm max} \leq \sz_{\rm max}$, which means the semigroup $P_t^f$ converges exponentially in the $L^2 (\nz)$-norm.

(ii). The next step, we show that the technical assumption (\ref{technical}) could be removed. Since $f\geq 1$, the exponential $f$-ergodicity of $P_t$ implies its exponential ergodicity:
$$
\lt\| P_t (i, \cdot) -\pz \rt\|_{\rm var} \leq C(i) e^{- \vz' t}, \qqd \forall t\geq 0, \  i \in E.
$$
Moreover, by \cite[Theorem 4.43]{Chen2004} and references therein, the constant $C(i)$ satisfies $C(i)\in L^1 (\pz)$. If $f$ is bounded, the ergodicity of $P_t$ ensures that the $f$-ergodicity holds. To be specific, let $f (i) \leq b$, $\forall i \in E$, then we have
\begin{align*}
\Big\| &  P_t (i, \cdot) - \pz \Big\|_{f} = \sup_{|g| \leq f} \lt| (P_t (i, \cdot)- \pz) g \rt| \leq \sup_{|g| \leq b} \lt| (P_t (i, \cdot)- \pz) g \rt| \\
&= \sup_{|g|/b \leq 1} b \lt| (P_t (i, \cdot) - \pz) (g/b) \rt| = b \| P_t (i, \cdot) -\pz \|_{\rm var},
\end{align*}
which means the exponential $f$-ergodicity holds with constant $C (i, f) = b C(i)$. Then, the conclusion holds by the method we used in the proof (i).

If the assumption (\ref{technical}) is invalid, we can define bounded functions as $f_N:= f\wedge N$, $N \in \mathds{N}^+$. Based on the above discussion, we have
$$
\lt\| \lt( P_t^{f_N} - \pz^{f_N} \rt) g \rt\|_{L^2 (\nz)}^2 \leq e^{-2 t \vz_{\rm max}}, \qqd \text{$\| g \|_{L^2 (\nz)} =1$}.
$$
Note that the right hand side of last inequality is independent of $N$. By dominated convergence theorem, we obtain the exponential convergence of the semigroup $P_t^f$ in the $L^2 (\nz)$-norm by letting $N\rightarrow \infty$.

(iii). Finally, we prove the necessity of Theorem \ref{Thm f-ergodic}. Assume $f\in L^2 (\pz)$ and $P_t$ is reversibility. If $P_t^f$ converges exponentially in the $L^2 (\nz)$-norm with $\sz_{\rm max} >0$, for any $0< s\leq t$, we have
\begin{align*}
\| P_s P_{t-s} &(i, \cdot) - \pz \|_f = \lt\| f \lt[ P_{t-s} \lt( \frac{P_s (i, \cdot)}{\pz_{\cdot}} -1  \rt) \rt] \rt\|_{L^1 (\pz)} \qqd \text{(by Lemma \ref{f-erg lmm4})} \\
&= \sum_{j \in E} \pz_j f(j) \lt| \sum_{k \in E} P_{t-s} (j, k) \lt( \frac{ P_{s} (i, k) }{\pz_k} -1 \rt) \rt| \\
&= \sum_{j \in E} \pz_j f^2 (j) \lt| \sum_{k \in E} \frac{1}{f(j)} P_{t-s} (j, k) f(k) \lt( \frac{1}{f(k)} \frac{ P_{s} (i, k) }{\pz_k } - \frac{1}{f(k)} \rt) \rt| \\
&= \sum_{j \in E}  \pz_j f^2 (j) \lt|P_{t-s}^f \lt( \frac{1}{f (\cdot)} \frac{ P_s (i, \cdot)}{ \pz_{\cdot} } - \frac{1}{f (\cdot)} \rt) (j) \rt| \\
&= \lt\| P_{t-s}^f \lt( h_s (i, \cdot) \rt) \rt\|_{L^1 (\nz)} \\
&\leq \lt\| P_{t-s}^f \lt( h_s (i, \cdot) \rt) \rt\|_{L^2 (\nz)} \pz \lt(f^2 \rt)^{1/2},
\end{align*}
The last step is Cauchy-Schwarz inequality and $h_s (i , \cdot)$ is defined as
\be\lb{h equ}
h_s (i, j) = \frac{1}{f(j)} \frac{ P_s (i, j)}{ \pz_j } - \frac{1}{f(j)}.
\de
For any $i\in E$ and $s>0$, we have $ \pz^f \lt( h_s (i, \cdot) \rt) =0$. By the exponential convergence of $P_t^f$ we have
$$
\lt\| P_{t-s}^f \lt( h_s (i, \cdot) \rt) \rt\|_{L^2 (\nz)} \leq e^{- \sz_{\rm max} (t-s)} \lt\|  h_s (i, \cdot) \rt\|_{L^2 (\nz)} ,
$$
where
\begin{align*}
\lt\| h_s (i, \cdot) \rt\|_{L^2 (\nz)}^2 &= \sum_{j \in E} \pz_j f^2 (j) \lt( \frac{1}{f (j)} \frac{ P_s (i, j)}{ \pz_j } - \frac{1}{f (j)} \rt)^2 \\
&= \sum_{j \in E} \lt( \frac{P_s (i, j)}{\pz_j} \rt)^2 \pz_j  - 1 \\
&= \frac{P_{2s} (i, i)}{ \pz_i } - 1,
\end{align*}
the last step depends on the reversibility of $P_t$. Hence, we obtain
\begin{align*}
\| P_t & (i, \cdot)  - \pz \|_f \leq \pz \lt( f^2 \rt)^{1/2} \lt\| P_{t-s}^f \lt( h_s (i, \cdot) \rt) \rt\|_{L^2 (\nz)} \\
&\leq \pz \lt( f^2 \rt)^{1/2} e^{- \sz_{\rm max} (t-s)} \lt\| h_s (i, \cdot) \rt\|_{L^2 (\nz)} \\
&= \pz \lt( f^2 \rt)^{1/2} e^{- \sz_{\rm max} t} \lt[ e^{\sz_{\rm max} s} \lt( \frac{P_{2s} (i, i)}{\pz_i} -1 \rt)^{1/2} \rt].
\end{align*}
Let $s \rightarrow 0$ and denote $C(i, f)$ by
$$
C(i, f) := \pz \lt( f^2 \rt)^{1/2} \lt( \frac{1}{\pz_i} -1 \rt)^{1/2},
$$
and then we get the exponential $f$-ergodicity
$$
\| P_t (i, \cdot)  - \pz \|_f \leq C(i, f) e^{- \sz_{\rm max} t}, \qqd \forall i \in E, \  t \geq 0,
$$
with $\sz_{\rm max} \leq \vz_{\rm max}$.
\deprf

Depending on Theorem \ref{Thm f-ergodic}, the only thing left to consider is the relationship between $\gap (Q)$ and $\sz_{\rm max}$. The crucial method is Poincar\'e inequality.

{\bf Proof of Theorem \ref{Thm f-ergodic1}.} For any $g$ satisfying $fg \in L^2 (\pz)$, we have $g \in L^2(\nz)$, and then the function
$$
F(t) = \lt\| \lt( P_t^f - \pz^f \rt) g \rt\|_{L^2 (\nz)}^2
$$
is well-define. Review the definition of the exponential convergence of $P_t^f$, we have $F(t) \leq F(0) e^{-2 \sz_{\rm max} t}$. Dividing by $t$, we get
\be\lb{thm =2}
\frac{\d}{\d t} F(t) \Bigg|_{t=0} \leq -2 \sz_{\rm max} F(0).
\de
By part (3) of Lemma \ref{f-erg lmm1}, we have
$$
F(0)= \lt\| g- \pz^f (g) \rt\|_{L^2 (\nz)}^2 = \pz \lt( f^2 g^2 \rt) - \pz^2 (fg)= \var_{\pz} (fg),
$$
and
$$
\frac{\d}{\d t} F(t) \Bigg|_{t=0} = \frac{\d}{\d t} \pz \lt(P_t^2(fg)\rt) \Bigg|_{t=0} = 2 (- Q (fg), fg)_{\pz}.
$$
Substituting these equations into (\ref{thm =2}), then
$$
\sz_{\rm max} \var_{\pz} (fg) \leq (- Q (fg), fg)_{\pz}, \qqd fg\in L^2 (\pz),
$$
which is Poincar\'e inequality. Since the spectral gap can be redefined as the optimal constant of the Poincar\'e inequality (cf. \cite[Chapter 9]{Chen2004}), then we have $\sz_{\rm max} \leq \gap (Q)$.

Conversely, assume $\gap(L) >0$. We use the same notations as aforesaid. Since $P_t (fg) \in L^2 (\pz)$, by Poincar\'e inequality, we have
$$
\gap (Q) \var_\pz \lt( P_t (fg) \rt) \leq \lt(- Q (P_t (fg)), P_t (fg) \rt)_{\pz},
$$
and then $2 \gap (Q) F(t) \leq - F' (t)$ for every $t \geq 0$. Using Gronwall lemma, we have $F(t) \leq e^{-2 \gap (Q) t} F (0)$. Therefore, $\gap(Q) \leq \sz_{\rm max}$.
\deprf

{\bf Proof of Proposition \ref{crl f-ergodic}.} Let $P_t^*$ be the dual semigroup of $P_t$, and its generator is denoted by $Q^*$. Using the $h$-transform, we can consider the convergence of semigroup $P_t^{* f}$. Similar to the Definition \ref{Def l2 con}, we denote the $L^2 (\nz)$-exponential convergence rate by $\sz_{\rm max}^*$. In the same way of the proof of Theorem \ref{Thm f-ergodic1}, we have $\sz_{\rm max}^* = \gap (Q^*) = \gap(Q)$. The second equality is base on \cite[Chapter 9]{Chen2004}.

It should be noted that Lemma \ref{f-erg lmm4} is still effective in the irreversible case. Assume $0< \gap(Q) = \sz_{\rm max}^*$, then the semigroup $P_t^{*f}$ converges exponentially in the $L^2 (\nz)$-norm. By the method in part (iii) of the proof of Theorem \ref{Thm f-ergodic}, we obtain
\begin{align*}
\| P_t & (i, \cdot) - \pz \|_f \leq \pz\lt(f^2\rt)^{1/2} \lt\| P_{t-s}^{* f} (h_s (i, \cdot)) \rt\|_{L^2 (\nz)} \\
& \leq \pz\lt(f^2\rt)^{1/2} e^{-\sz_{\rm max}^* t} \lt[ e^{\sz_{\rm max}^* s} \lt(\frac{P_{2s} (i, i) }{\pz_i} -1 \rt)^{1/2} \rt],
\end{align*}
where $h_s (x, \cdot)$ is defined as (\ref{h equ}). Let $s \rightarrow 0$, and then we get the $f$-ergodicity immediately, which satisfies $\gap(Q) \leq \vz_{\rm max}$.
\deprf

\noindent {\bf Acknowledgements} Research supported in part by the National Natural Science Foundation of China (No. 11701588, 61773411).

\end{document}